\documentclass[3p,times]{elsarticle}

\usepackage{amssymb}

\usepackage[figuresright]{rotating}

 \newtheorem{theorem}{Theorem}[section]
\newtheorem{lemma}[theorem]{Lemma}
\newtheorem{e-proposition}[theorem]{Proposition}

\newtheorem{e-definition}[theorem]{Definition\rm}

\newtheorem{example}{\bf Example\/}

\begin{document}

\begin{frontmatter}




\title{On a problem of J. Nakagawa, K. Sakamoto, M. Yamamoto}
\author{HNAIEN Dorsaf }
\ead{hnaien.dorsaf@gmail.com}
\author{KELLIL Ferdaous }
\ead{kellilferdaous@yahoo.fr,\,ferdaous.kellil@fsm.rnu.tn}
\author{LASSOUED Rafika}
\ead{rafika.lassoued@gmail.com}

  \address{Laboratory of Applied Mathematics and Harmonic Analysis, 3038 Sfax, Tunisia}




\begin{abstract}
In this paper, we give a positive answer to a problem posed by Nakagawa, Sakamoto and Yamamoto concerning a nonlinear equation with a fractional derivative.
\end{abstract}


\end{frontmatter}


\section{Introduction}
\par\indent In their overview paper concerning the mathematical analysis of fractional equations, Nakagawa, Sakamoto and Yamamoto \cite{1} posed the problem concerning global solutions and blowing-up in a finite time of solutions to the equation
\begin{equation}\label{6}
\left\{ \begin{array}{cll}
  ^{C}D_{0_{+}}^{\alpha} u(t)& = & -u(t)(1-u(t)),~~t>0,\\
  u(0) & = & u_{0},
\end{array}\right.
\end{equation}
where $^{C}D_{0_{+}}^{\alpha}$ is the Caputo derivative defined for $g \in C^{1}[0,T]$ by
$$ ^{C}D_{0_{+}}^{\alpha}g(t)=\frac{1}{\Gamma(1-\alpha)} \int_{0}^{t} (t-\tau)^{-\alpha} g'(\tau) \,d\tau ,$$
for $0<\alpha <1$.\\
\par\indent Let us recall, in the case $\alpha =1$, the results concerning solutions of $(\ref{6})$ :\\ \\
- For $0< u(0)< 1$, the solution exists globally. Moreover, $$|u(t) | \leq \frac{1}{e^{t}(1-u_{0})} \longrightarrow 0, \hbox{ as }{t  \longrightarrow + \infty}. $$\\
- For $u(0) > 1$, the solution can not exist globally.\\ \\
Here, we show that the same conclusions are valid for equation $(\ref{6})$. Moreover we analyse : \\
\begin{enumerate}
\item The large time behavior of the global solution.\\
\item The blow-up time and profile of the blowing-up solutions. \\
\end{enumerate}

\par\indent Note that if we set $w=u-1$, then $(\ref{6})$ reads
$$^{C}D_{0_{+}}^{\alpha}w(t)  =  w(t)(1+w(t)),$$
which describes the evolution of a certain species; the reaction term $w(1+w)$ describes the law of increase of the species.

\section{ Preliminaries }
\label{}
\par\indent In this section, we present some definitions and results concerning fractional calculus that will be used in the sequel. For more information see \cite{2}.
\par\indent The Riemann-Liouville fractional integral of order $ 0 < \alpha <1$  of the integrable function $f : \mathbb{R}^{+} \longrightarrow  \mathbb{R}$ is
$$  J_{0_{+}}^{\alpha} f(t):= \frac{1}{\Gamma(\alpha)} \int_{0}^{t} (t-\tau)^{\alpha-1} f(\tau) \,d\tau, ~~~t>0,$$
where $\Gamma (\alpha)$ is the Euler Gamma function.
\par\indent The Riemann-Liouville fractional derivative of an absolutely continous function $f(t)$ of order $0 <\alpha <1$ is
$$ D_{0_{+}}^{\alpha} f(t):= \frac{1}{\Gamma(1- \alpha)} \frac{d}{dt} \int_{0}^{t} (t-\tau)^{- \alpha} f(\tau) \,d\tau .$$
\par\indent The Caputo fractional derivative of an absolutely continous function $f(t)$ of order $0 <\alpha <1$ is defined by
$$ ^{C}D_{0_{+}}^{\alpha} f(t):=  J_{0_{+}}^{1-\alpha} f'(t)= \frac{1}{\Gamma(1-\alpha)} \int_{0}^{t} (t-\tau)^{-\alpha} f'(\tau) \,d\tau .$$
\par\indent Both derivatives present a drawback : \\ \\
- The Riemann-Liouville derivative of a constant is different from zero,
$$ D_{0_{+}}^{\alpha} C \neq 0 ,$$
while the Caputo derivative require $f'(t)$ to calculate $^{C}D_{0_{+}}^{\alpha} f(t)$, for $0 < \alpha <1$.\\ \\
- We know that the Riemann-Liouville derivative of the Weierstrass function exists for any $0< \alpha < 1$, but not for $\alpha = 1 $.\\ \\ But for regular function with $f(0)= 0$, both definitions coincide.\\
\par\indent Next, we recall a lemma that will be used hereafter.
\begin{lemma}( see \cite{3}). Let $a,~ b,~ K,~\psi $ be non-negative continuous functions on the interval $I= ( 0, T )~~( 0 < T \leq \infty)$, let $ \omega : (0, \infty) \longrightarrow \mathbb{R}$ be a continuous, non-negative and non-decreasing function with $\omega(0)= 0$ and $\omega(u)>0$ for $u >0$, and let $A(t)= \max_{0 \leq s \leq t} a(s)$ and $B(t)=\max_{0 \leq s \leq t} b(s)$. Assume that
$$ \psi(t) \leq a(t)+b(t) \int_{0}^{t} K(s) \omega( \psi(s))\, ds,~~ t \in I. $$
Then
$$  \psi(t) \leq H^{-1} \Big[ H(A(t))+ B(t) \int_{0}^{t} K(s) ds   \Big],~~ t \in (0, T_{1}),  $$
where $H(v)= \int_{v_{0}}^{v}  \frac{d\tau}{\omega(\tau)}~~(v \geq v_{0} >0)$, $H^{-1}$ is the inverse of $H$ and $T_{1} >0$ is such that $H(A(t))+ B(t) \int_{0}^{t} K(s) ds  \in D(H^{-1})$ for all $t \in (0, T_{1})$.
\end{lemma}
\par\indent Here, we consider the problem
\begin{equation}\label{P}
\left\{ \begin{array}{cll}
  ^{C}D_{0_{+}}^{\alpha} u(t)& = & -u(t)(1-u(t)),\\
  u(0) & = & u_{0},
\end{array}\right.
\end{equation}
for $0 <\alpha <1$ and $u_{0}>0$.
\section{Main results}
\par\indent The local existence of solutions to $(\ref{P})$ is assured by the
\begin{theorem}(see \cite{2}).
 We consider the fractional differential equation of Caputo's type given by
\begin{equation}\label{1}
\left\{ \begin{array}{cll}
  ^{C}D_{0_{+}}^{\alpha} u(t)& = & f(t,u(t)),~~t>0,\\
  u(0) & = & u_{0}.
\end{array}\right.
\end{equation}
For $0<\alpha<1$, $u_{0}\in\mathbb{R},~~b>0$ and $T>0$. \\
Assume that
\begin{enumerate}
\item $f\in C(R_{0},\mathbb{R})$ where $R_{0}=\{(t,u),~~0\leq t\leq T,~~|u-u_{0}|\leq b\}$ and $|f(t,u)|\leq M $ on $R_{0}$;
\item $|f(t,u)-f(t,v)|\leq L|u-v|,~~L>0,~~(t,u)\in R_{0}$.
\end{enumerate}
Then there exists a unique solution $u \in C([0,h])$ for $(\ref{1})$, where $h=\min \Big\{T,\displaystyle \Big(\frac{b\Gamma(\alpha+1)}{M}\Big)^{\frac{1}{\alpha}} \Big\}$.
\end{theorem}
\begin{theorem}
Let $u$ be the solution of problem $(\ref{P})$. We have :\\ \\
- If $~0 < u_{0} < 1$, the solution is global and it satisfies $~0 < u <1$. Moreover, $u$ is given by
$$ u(t)=E_{\alpha}(-t^{\alpha}) u_{0}+\int_{0}^{t}(t-s)^{\alpha-1}E_{\alpha,\alpha}(-(t-s)^{\alpha})~u^{2}(s)ds,$$
and for some constants $c>0$ and $c_{1}>0$, we have
$$ 0 < u(t) \leq \frac{1}{\frac{1}{c u_{0}}- \frac{c_{1}}{\alpha} t^{\alpha}},~~~~ t > T_{0}:= \Big( \frac{\alpha}{c_{1} c u_{0}} \Big)^{\frac{1}{\alpha}}.$$
- If $~u_{0}> 1$, the solution blows-up in a finite time $T^{*}$ : $ \displaystyle \lim_{t \longrightarrow T^{*}} u(t)= + \infty$. \\
Moreover, we have the bilateral estimate :
$$\overline{ w }(t)+1 \leq u(t) \leq \widetilde{w}(t)+1,    $$
and
$$\Big( \frac{  \Gamma (\alpha+1)}{ 4 (u_{0}-\frac{1}{2}) } \Big)^{\frac{1}{\alpha}}  \leq T^{*} \leq \Big( \frac{  \Gamma (\alpha+1)}{ u_{0}-1} \Big)^{\frac{1}{\alpha}} ,$$
where

$$ \widetilde{w}(t)+ \frac{1}{2} \sim \frac{\Gamma(2\alpha)}{\Gamma(\alpha)} (T_{\widetilde{w}} - t )^{- \alpha}, \hbox{~~as}~~ t \longrightarrow  T_{\widetilde{w}},$$

$$ \overline{w}(t) \sim \frac{\Gamma(2\alpha)}{\Gamma(\alpha)} (T_{\overline{w}} - t )^{- \alpha}, \hbox{~~as}~~ t \longrightarrow  T_{\overline{w}}.$$
Here, $T_{\widetilde{w}}$ is the blow-up time of $~\widetilde{w}$, which satisfies
$$\Big( \frac{  \Gamma (\alpha+1)}{ 4 (u_{0}-\frac{1}{2})} \Big)^{\frac{1}{\alpha}} \leq T_{\widetilde{w}} \leq \Big( \frac{  \Gamma (\alpha+1)}{u_{0}-\frac{1}{2}}\Big)^{\frac{1}{\alpha}},$$
and $T_{\overline{w}}$ is the blow-up time of  $~\overline{w}$, which satisfies
$$\Big( \frac{  \Gamma (\alpha+1)}{ 4 (u_{0}-1)} \Big)^{\frac{1}{\alpha}} \leq T_{\overline{w}} \leq \Big( \frac{  \Gamma (\alpha+1)}{u_{0}-1}\Big)^{\frac{1}{\alpha}}.$$

\end{theorem}

\textbf{Proof of Theorem 3.2.} \\
\par\noindent \textbf{Part 1}. If  $0 < u_{0} < 1$, then the solution is global.\\
The solution to $(\ref{P})$ is given by
\begin{equation}\label{3}
u(t)=E_{\alpha}(-t^{\alpha}) u_{0}+\int_{0}^{t}(t-s)^{\alpha-1}E_{\alpha,\alpha}(-(t-s)^{\alpha})~u^{2}(s)ds.
\end{equation}
Where the Mittag-Leffler functions  $E_{\alpha}(-t^{\alpha})$  and  $E_{\alpha,\alpha}(-t^{\alpha})$  are defined by :
$$E_{\alpha}(-t^{\alpha}) = \sum_{j=0}^{\infty} \frac{(-1)^{j} t^{\alpha j}}{ \Gamma(\alpha j +1)},$$
$$E_{\alpha,\alpha}(-t^{\alpha}) = \sum_{j=0}^{\infty} \frac{(-1)^{j} t^{\alpha j}}{ \Gamma(\alpha j +\alpha)}.$$
If $ u_{0}> 0$, then $ u(t) >0$ as $E_{\alpha}(-t^{\alpha}) > 0$ and $E_{\alpha,\alpha}(-t^{\alpha}) > 0$.\\ \\
Now, we set the function $\overline{u}(t)=1, ~~t>0$.\\ \\
As $ 0 < u_{0} < 1$, then $u_{0} < \overline{u}(0).$ In addition, we have
$$^{C}D_{0_{+}}^{\alpha}\overline{u}(t)=0= -\overline{u}(t)(1-\overline{u}(t)).$$
Hence $\overline{u}$ is an upper solution of the equation $(\ref{P})$, and we have $ u(t) < \overline{u}(t)=1$, (see \cite{4}, Thm. 2.4.3, p. 32). \\
\par\indent Now, we examine the large time behavior of the global solution $0< u < 1$.\\
For, let us recall the estimates ( see \cite{5}) :\\
- For $0< \alpha <1$, there exists a constant $c >0$ such that,
\begin{equation}\label{9}
0< E_{\alpha}(-t^{\alpha})  \leq \frac{c}{1+ t^{\alpha}} \leq c,~~~~t >0.
 \end{equation}
- For $0< \alpha <1$, there exists a constant $c_{1} >0$ such that
\begin{equation}\label{10}
0<  t^{\alpha-1}E_{\alpha,\alpha}(-t^{\alpha}) \leq  c_{1}t^{\alpha-1},~~~~t >0.
 \end{equation}
From $(\ref{3})$ and using the inequalities $(\ref{9})$ and $(\ref{10})$, we obtain
\begin{equation}\label{11}
 u(t) \leq  c u_{0}+ c_{1} \int_{0}^{t}(t-s)^{\alpha-1}  u^{2}(s)ds  .\end{equation}
We apply Lemma 2.1 to $( \ref{11})$ with $\omega(x)=x^{2},~K(s)=(t-s)^{\alpha-1},~ A(t)=c u_{0},~ B(t)= c_{1}.$\\
For $t > T_{0}$, we have
$$ H(c u_{0})+ \frac{ c_{1}}{\alpha} t^{\alpha} \in D(H^{-1}),$$
where $\displaystyle H(v)= \frac{1}{v_{0}}-\frac{1}{v}~$ and $~\displaystyle H^{-1}(z)= \frac{1}{\frac{1}{v_{0}} -z},~~ z\neq \frac{1}{v_{0}}$.\\
So we obtain,
$$ u(t) \leq  H^{-1}\Big[H(c u_{0})+ \frac{ c_{1}}{\alpha} t^{\alpha}\Big].$$
Therefore
$$u(t) \leq \frac{1}{\frac{1}{c u_{0}}- \frac{c_{1}}{\alpha} t^{\alpha}},~~ t > T_{0}.$$

\par\noindent \textbf{Part 2.} If $ u_{0} > 1$, then the solution blows-up in a finite time.\\
\begin{enumerate}
\item We show that $ u > 1$. For, let us define the new unknown function $w = u-1$. The function $w$ satisfies
  \begin{equation}\label{4}
\left\{ \begin{array}{cll}
  ^{C}D_{0_{+}}^{\alpha}w(t) & = & w(t)(1+w(t)),\\
  w(0):=w_{0} & = & u_{0}-1.
\end{array}\right.
\end{equation}
As $u_{0}>1$, then $ w_{0}>0$.
Moreover, we have (\cite{2})
$$w(t)=E_{\alpha}(t^{\alpha}) w_{0} +\int_{0}^{t}(t-s)^{\alpha-1}E_{\alpha,\alpha}((t-s)^{\alpha})~w^{2}(s)ds.$$
Therefore, $w>0$; hence $u>1$.\\
\item We prove that $ u$ blows-up in a finite time.\\ \\
Since we have $w(t)= u(t)-1$, it is seen that if $u(t) \longrightarrow \infty$ as $ t \longrightarrow T^{*}$, then $w(t)\longrightarrow \infty$ as $ t \longrightarrow T^{*} $ and vice versa. That is $w$ and $u$ will have the same blow-up time.\\
We now must examine the blow-up properties of $w$, the solution of problem $(\ref{4})$. These are obtained by comparing $w(t)$ with the solutions of the following problems :
 \begin{equation}\label{7}
\left\{ \begin{array}{cll}
  ^{C}D_{0_{+}}^{\alpha}\overline{w}(t) & = & \overline{w}~^{2}(t),\\
  \overline {w}(0) & = & w_{0},
\end{array}\right.
\end{equation}

and

\begin{equation}\label{8}
\left\{ \begin{array}{cll}
  ^{C}D_{0_{+}}^{\alpha}\widetilde{w}(t) & = & ( \widetilde{w}(t)+ \frac{1}{2})^{2},\\
  \widetilde{w}(0) & = & w_{0}.
\end{array}\right.
\end{equation}
We see by comparaison  (\cite{4}) that $$ \overline{w}(t) \leq w(t) \leq  \widetilde{w}(t),~~ 0\leq t < \min \{ T_{\overline{w}},T_{\widetilde{w}} \}.$$
Following the paper of Kirk, Olmstead and Roberts \cite{6}, we may assert that the solution $\overline{w}$ ( resp. $\widetilde{w}$) blows-up in a finite time
$T_{\overline{w}}$ ( resp. $T_{\widetilde{w}}$), such that
$$\Big( \frac{  \Gamma (\alpha+1)}{ 4 w_{0}} \Big)^{\frac{1}{\alpha}} \leq T_{\overline{w}} \leq \Big( \frac{  \Gamma (\alpha+1)}{w_{0}}\Big)^{\frac{1}{\alpha}},$$
and
$$\Big( \frac{  \Gamma (\alpha+1)}{ 4 (w_{0}+\frac{1}{2})} \Big)^{\frac{1}{\alpha}} \leq T_{\widetilde{w}} \leq \Big( \frac{  \Gamma (\alpha+1)}{w_{0}+\frac{1}{2}}\Big)^{\frac{1}{\alpha}}.$$
So we have the following estimates
$$T_{\widetilde{w}} \leq T^{*} \leq T_{\overline{w}}.$$
Whereupon
$$\Big( \frac{  \Gamma (\alpha+1)}{ 4 (w_{0}+\frac{1}{2})} \Big)^{\frac{1}{\alpha}}  \leq T^{*} \leq \Big( \frac{  \Gamma (\alpha+1)}{ w_{0}} \Big)^{\frac{1}{\alpha}} .$$
$\hfill \Box$
\end{enumerate}
\section{ Numerical implementation}
In this section, we will approximate the solution $u$ given by $(\ref{3})$. For, we need a numerical approximation of the convolution integral; this can be obtained using the convolution quadrature method.\\
\par As it has been explained in \cite{7}, a convolution quadrature approximates the continuous convolution
$$  \int_{0}^{t} K(t-s) f(s) \, ds, ~~ t >0,$$
by a discrete convolution with a step size $h >0$.
Then $$\int_{0}^{t_{n}} K(t_{n}-s) f(s) \, ds \sim  \sum_{j=0}^{n} \omega_{n-j} f(t_{j}), $$
where $t_{j}= j h, ~ j=0,1,2,...,n$ and the convolution quadrature weights $\omega_{j}$ are determined from their generating power series as
$$ \sum_{j=0}^{\infty} \omega_{j} \zeta^{j} = \mathcal{L} \Big\{ K(t): \frac{\delta(\zeta)}{h}  \Big\} .  $$
Here $\mathcal{L} \{K(t):s \}$ is the Laplace transform of $K(t)$ and $\delta (\zeta)$ is the generating polynomial for a linear multistep method.
 \par Let $u_{n} $ be the approximation of $u(t_{n})$ for $n \geq 0$. Using the convolution quadrature method we obtain
$$  u_{n} = (1- \omega_{0})^{-1} \Big[ E_{ \alpha}(- t^{ \alpha}) u_{0}+ \sum_{j=0}^{n-1} \omega_{n-j} u_{j} \Big],~~n=1,2,3....$$

\begin{figure}
\begin{center}
  \includegraphics[width=10cm]{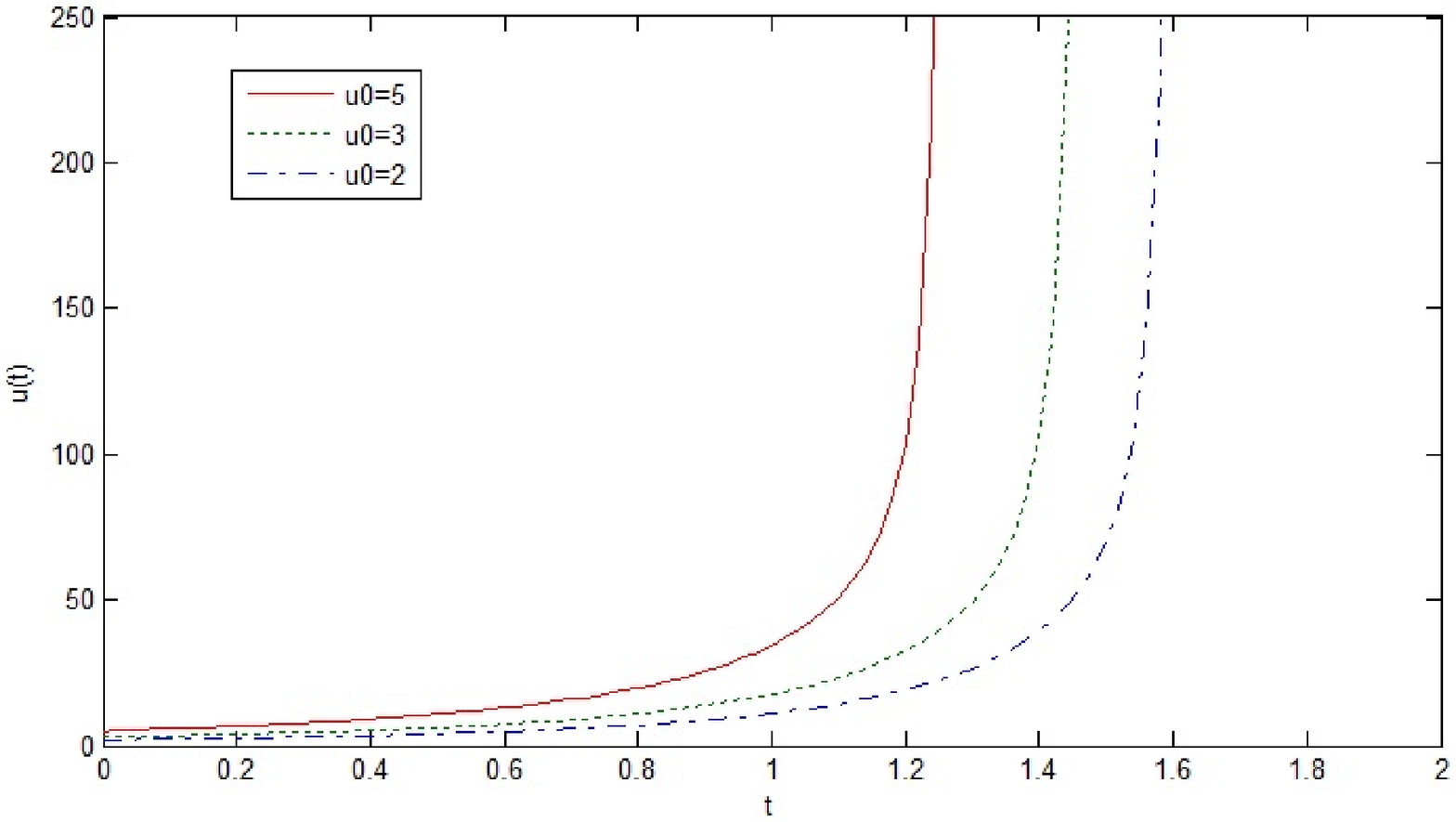}\\
  \caption{solutions for $\alpha = 0.5$ and $u_{0}= 5,~3,~2$. } \label{exp un}
  \includegraphics[width=10cm]{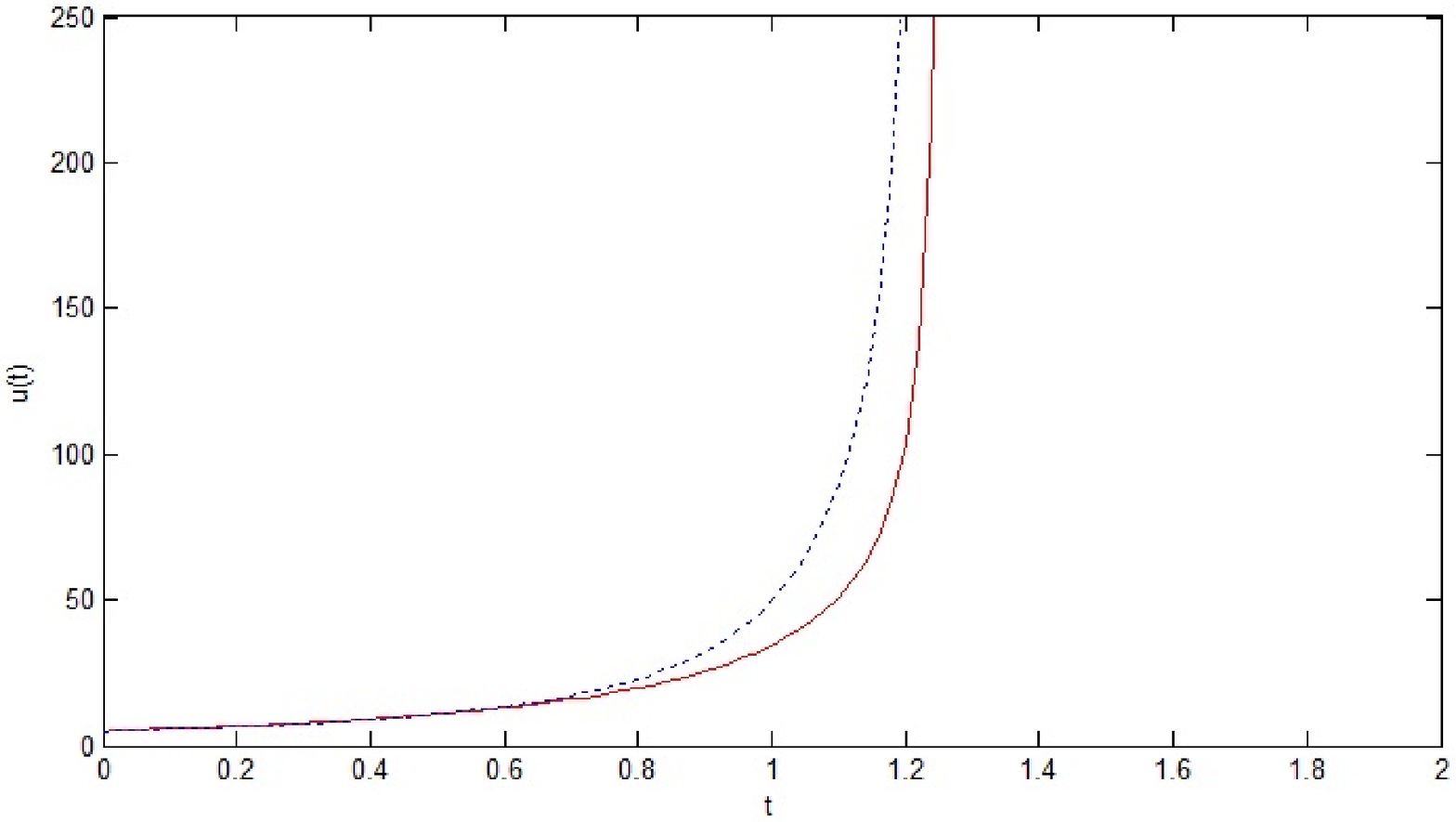}\\
  \caption{solutions for $u_{0}= 5$ and $\alpha = 0.3,~0.5$. } \label{exp deux}
 \end{center}
\end{figure}

\par Now, we introduce the following algorithm which gives the numerical approximation of solution to equation $(\ref{P})$.\\ \\
\textbf{ Algorithm }\\
\textbf{Input : } Give $\alpha$, $0 < \alpha <1$ and $u_{0}, ~~ u_{0} >1$.\\
\textbf{ Initializations :} Discretize the time with a step size $h >0$; $t_{i}= ih$, for all $ i= 1,2,...,n$, $u^{1}_{appx}= u_{0}, ~~ u^{1}= (u_{0})^{2}. $   \\
\textbf{Step 1 :} Approximate the Mittag-Leffler function \textbf{GML}.\\
\textbf{Step 2 :} Calculate convolution quadrature weights \textbf{W} using the fast Fourier transform (FFT).\\
\textbf{Step 3 :} Calculate $ u^{i}_{appx}$.\\
\textbf{do}

$\begin{array}{l}
  u^{i}= \textbf{GML} * u^{1}_{appx} + \textbf{W} * u^{i-1} . \\
  u^{i}_{appx} = (1- \textbf{W}(1) )^{-1} * u^{i}. \\
  u^{i }= ( u^{i}_{appx})^{2}.\\
   i= i+1.
\end{array}$\\
\textbf{until} ($u^{i}_{appx}$ blows up) or ($ i > n$ ).\\
\textbf{ Output :} Numerical approximation of $u$.\\
\begin{example}
For $Figure 1$, we set $\alpha=0.5$; the initial conditions are respectively $u_{0}=5,~~u_{0}=3$ and $u_{0}=2$.\\
For $Figure 2$, we take the initial condition $u_{0}=5$ and we plot the solutions; the dotted curve is the solution for $\alpha=0.3$ and the solid curve corresponds to the solution for $\alpha=0.5$.\\
\par As it has been proved, the solution blows up in a finite time which depends on $u_{0}$ and $\alpha$.

\end{example}





\bibliographystyle{elsarticle-num}



\end{document}